\documentclass[12pt,reqno]{amsart}
\usepackage[T1]{fontenc}
\usepackage[utf8]{inputenc}

\usepackage{comment} 
\usepackage{xcolor}
\usepackage{amssymb}
\usepackage{graphicx}
\usepackage{subfigure}
\usepackage{float} 
\graphicspath{{fig/}}
\usepackage[margin=3cm]{geometry}

\usepackage{tikz-cd}

\def\sloppy{\tolerance500 \emergencystretch 3em\relax}

\let\ts=\textstyle
\let\ol=\overline
\let\ep=\varepsilon

\let\del=\partial
\let\joli=\mathcal
\let\ds=\displaystyle
\let\ts=\textstyle
\def\delz{{\del\over\del z}}
\def\delt{{\del\over\del t}}
\def\delZ{{\del\over\del Z}}
\def\zbar{{\ol z}}
\def\Cc{\mathbb C}
\def\Dd{\mathbb D}
\def\R{\mathbb R}
\def\N{\mathbb N}
\def\Hh{\mathbb H}
\def\zbar{{\overline z}}
\def\inv{^{-1}}

\catcode`@=11
\def\paragraph{\@startsection {paragraph}{4}\z@{5pt} {-\fontdimen 2\font }\itshape}
\catcode`@=12

\newtheorem{theo}{Theorem}[section]
\newtheorem{cor}[theo]{Corollary}
\newtheorem{lemma}[theo]{Lemma}
\newtheorem{prop}[theo]{Proposition}
\theoremstyle{definition}
\newtheorem{deff}[theo]{Definition}
\newtheorem{rem}[theo]{Remark}

\newtheorem{notation}[theo]{Notation}

\theoremstyle{definition}

\newcommand{\ord}[0]{\text{ord}\,}

\title{Generic Antiholomorphic Polymonial Vector Fields}
\author[J. Godin]{Jonathan Godin} 
\address{Universit\'e de Moncton, Campus de Shippagan, 218 
    boul.~J.-D.-Gauthier, Shippagan (NB), E8S 1P6, Canada}
\email{jonathan\_godin@msn.com}
\author[J. Perazzelli]{Jérémy Perazzelli}
\address{D\'epartement de math\'ematiques et de statistique, Universit\'e 
    de Montr\'eal, C.P. 6128, Succursale Centre-ville, Montr\'eal (QC), 
    H3C 3J7, Canada}
\email{jeremy.perazzelli@umontreal.ca} 
\date{\today}

\usepackage[hidelinks]{hyperref}
\hypersetup{
    colorlinks,
    linkcolor={red!50!black},
    citecolor={blue!50!black},
    urlcolor={blue!80!black}
}

2

\begin{document}
\expandafter\let\expandafter\over\csname @@over\endcsname
\def\d{\,{\rm d}}

\begin{abstract}
  An analytic classification of generic anti-polynomial vector fields $\dot z = \overline{P(z)}$ is given in term of a topological and an analytic invariants. The number of generic strata in the parameter space is counted for each degree of $P$. A Realization Theorem is established for each pair of topological and analytic invariants. The non-generic case of a maximal number of heteroclinic connections is also given a classification. The bifurcation diagram for the quadratic case is presented.
\end{abstract}
\maketitle

\tableofcontents

\section{Introduction}
In this paper, we are mainly concerned with the study and classification
of vector fields $\ol{P(z)}\delz$, where $P$ is a polynomial and
$\ol{\,\cdot\,}$ is the complex conjugation. By classifying, we mean 
describing the equivalence classes of the family of vector fields
$\ol{P(z)}\delz$ under some equivalence relation, such as the topological
orbital equivalence or holomorphic equivalence.

This article follows in the step of the visionary work of
Douady, Estrada and Sentenac in their unpublished paper~\cite{DES} on
structurally stable polynomial vector fields, and later the work
of Branner and Dias~\cite{BD}, expanding the analysis to every polynomial
vector fields. The classification of real-coefficients polynomial vector
fields has recently been studied in~\cite{GR}.

Early in the development of complex analysis, anti-holomorphic vector fields
found many applications in physics because of complex potentials; broadly
speaking, given a (possibly ramified) holomorphic function $f\colon U\to \Cc$
(the complex potential), over some open set $U \subseteq \Cc$, the dynamics of
$\ol{f'(z)}\delz$ may be thought of as the laminar flow of a Newtonian fluid in
the plane. Although the flow of a real fluid is much more complicated than
this, these simple potential flows can be used to construct a variety
of flows, by superposition, patching or otherwise, that are much closer
to what is observed in the physical world. Potential flows have been
used in the early development of aerodynamics~\cite{van} and in recent
work (e.g.~\cite{CKM}). From this point of view, generic anti-holomorphic 
polynomials arise as vector fields on $\Cc$ that are 
\begin{itemize}
    \item meromorphic on the whole Riemann sphere;
    \item and have only simple singularities (forced to be saddle points) in
    $\Cc$ and no homoclinic loops or heteroclinic connections.
\end{itemize}
Such a flow with $k+1$ saddle points is characterized at the topological level by 
the organization of $k$ strips (zone isomorphic to horizontal strip) in which the 
fluid flows from infinity to infinity. 
In Section~\ref{sec:combinatoire}, we prove that there are $A^r(k+2)$ such configurations
(see Equation~\eqref{eq:6}),
where $A^r(n)$ is the number of non crossing trees with $n$ vertices modulo
rotations of order $n$. Then the 
circulation and the flux of each strip determines completely the polynomial.

There has also been recent development in approximation of potential flows in
planar domains~\cite{lightning}. Rational functions are used to approximate the
flow. The method uses a cluster of poles near the corner to achieve very fast
approximation when the number of poles increases exponentially. Recent
work~\cite{tref} on the approximation of solutions of the Laplace equation
$\Delta f = 0$ in planar domains showed that rational approximations perform
better than polynomial approximations on some non-convex regions.

The most degenerate singularity is often the organizing centre of the dynamics.
The unfoldings of a multiple saddle point can be analyzed through the use of
a family of anti-polynomials (i.e.~polynomial in $\zbar$) or, equivalently, a family of rational function
of the form ${1\over P}$. The goal of the present paper is to classify
families of generic anti-polynomials.

On a different note, we also see this work as part of the analysis
of (holomorphic) polynomial vector fields on $\Cc$. Indeed, it is well known
that such a vector field has a pole at infinity, which is the organizing centre
of the dynamics. To study the unfoldings of this singularity, it is natural
to consider ${1\over P(z)}\delz$, which leads to $\ol{P(z)}\delz$ by observing
that
$$
    \dot z = {1\over P(z)} = {\ol{P(z)}\over |P(z)|^2}.
$$
In other words, the dynamics of ${1\over P(z)}\delz$ and $\ol{P(z)}\delz$
have the same orbits, but differ by a time
rescaling. In fact, in this paper, we will often consider both $\ol{P(z)}\delz$
and ${1\over P(z)}\delz$, and use whichever one is more convenient given
the context.

A different point of view can also be used to see this work as part of holomorphic
dynamics. With the idea that $w = \zbar$, we can see the ODE $\dot z  = \ol{P(z)}$
as a system of two equations in $\Cc^2$
$$
    \begin{cases}
        \dot z = \ol{P(\ol w)}, \\
        \dot w = P(z).
    \end{cases}
$$ 
This family of systems is invariant under the symmetry $(z,w)\mapsto (\ol w,\ol z)$.
The dynamics studied in this paper corresponds to the dynamics for real time in the 
invariant plane $\{w = \zbar\}$ in $\Cc^2$. See~\cite{ilyashenko} for more on
analytic systems.

The paper is organized as follow. In Section \ref{sec:desc}, we study the separatrix graph
and the two types of components in the complement of this graph. We also discuss the topological properties of the phase portrait. In Section~\ref{sec:combinatoire}, we count the number of equivalence classes under topological
orbital equivalence. Section~\ref{sec:analytique} is devoted to the analytic invariants of the
vector fields. In Section~\ref{sec:realisation}, we prove a realization theorem, proving that
every phase portrait is completely determined by the topological and
analytic invariants. Lastly, in Section~\ref{sec:non stable}, we initiate the study of
non-structurally stable vector fields.

\section{Description of the phase portrait}\label{sec:desc}
A polynomial $P(z) = a_{k+1} z^{k+1} + \cdots + a_0$ is monic if
$a_{k+1} =1$ and centred if $a_k = 0$.

\begin{notation}
Through the paper, we will denote by $\Xi_{k+1}$ the set of monic, centred
anti-polynomial vector fields of degree $k+1$ and we will implicitly identify any such vector
field with its defining polynomial. We will note such a vector field $\ol P \delz$, with
\begin{align*}
    \ol{P(z)} = \ol{z^{k+1} + a_{k-1}z^{k-1} + \cdots + a_0}.
\end{align*}
\end{notation}

For $\ol{P}\in\Xi_{k+1}$, we consider the system
\begin{equation}\label{eq:1}
    \dot{z} = \overline{P(z)}.
\end{equation}
Every singular point must a saddle point; this can be checked directly in
polar coordinates in a neighbourhood of the singular point. Therefore, 
no nodes, parabolic points or centres are possible in $\Cc$. Homoclinic loops or
limit cycles cannot occur, since they must have a positive index number. For the
same reason, a domain cannot be bounded by orbits.

\begin{figure}[h!]
    \subfigure[Outside a large disk in $\Cc$]{\includegraphics[scale=.8]{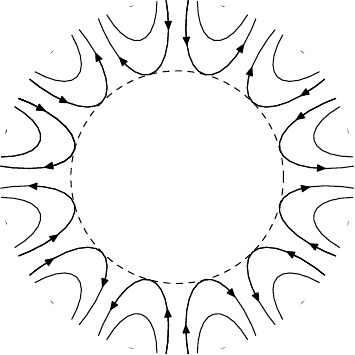}}
    \hskip2cm
    \subfigure[Near the circle at infinity in~$\ol \Cc$]{\includegraphics[scale=.8]{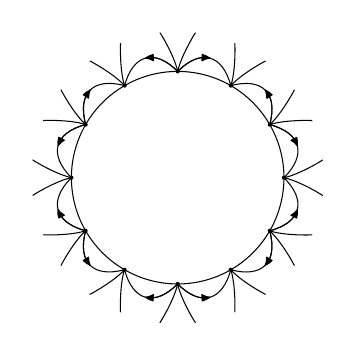}}
    \caption{Phase portrait of $\dot z = \ol{P(z)}$, with $\deg P = 5$}
    \label{fig:portrait loin}
\end{figure}

Near infinity, the behaviour of the dynamics is of the form $\dot w =
-\frac{w^{k+3}}{|w|^{2(k+1)}}(1+o(1))$, with $w = \frac 1z$ and $|w|$ small.
Hence, the system $\dot{w}$ has a parabolic point of codimension $k+2$ near the
origin. Consider the blow-up $\ol \Cc = \Cc \sqcup S^1$ at $\infty$ of the
Riemann sphere $\hat{\Cc}$ with fiber $S^1$ (one point for every oriented
direction at infinity). In this space, we can make sense of the notion of going
to $\infty$ in the direction $\theta$ for every $\theta\in S^1$. The blown up
parabolic point at infinity has $k+2$ attracting points and $k+2$ repulsing
points on $S^1$ in alternating order. These marked points are located at
$\exp\!\left( {2ij\pi \over 2k+4}\right)$ on the circle at infinity, with an
attractive (resp.~repulsive) point for $j$ even (resp.~$j$ odd). See
Figure~\ref{fig:portrait loin}.

\begin{deff}\label{def:eq-top-orb}
    Two vector fields $\ol{P_1}\delz$ and $\ol{P_2}\delz$ are said to be
    topologically orbitally equivalent if there exists a homeomorphism
    $h:\Cc\to\Cc$ that maps orbits of $\ol{P_1}\delz$ into orbits of
    $\ol{P_2}\delz$. In this paper, we further require that the homeomorphism
    fixes the order of the marked points at infinity (see above).
\end{deff}

The phase portrait may contain a finite number of simple saddle points,
multiple saddle points and heteroclinic connections. The latter two are
the only possible non-structurally stable objects. 


\subsection{Separatrix graph} The separatrix graph, formed of singular points
and separatrices, was already used in~\cite{DES} in the case of holomorphic
polynomial vector fields, so it is no surprise that this object completely
determines the topological structure of the phase portrait in the case at hand.
At the end of Section~\ref{sec:desc}, we will reduce the study of the structurally
stable phase portrait under topological orbital equivalence to that of
different graphs, deduced from the separatrix graph, called noncrossing trees.
The discrete and combinatorial properties of those trees were already studied
in \cite{NOY98}.

The basin of attraction of $\infty$ consists of a finite number of connected and simply connected regions of $\Cc$ bounded by singular points and orbits with bounded $\omega$-limits. If we know the basin of attraction at $\infty$ and its boundary, we can reconstruct the whole phase portrait. This allows us to reduce our problem to the study of the separatrix graph (defined in the following paragraph).


\begin{deff}\begin{enumerate}
    \item An incoming (resp.~outgoing) separatrix is an orbit of a non singular point
        on the boundary of the basin of attraction (resp.~repulsion) of the
        point at infinity. 
    \item We will call \emph{separatrix graph} of $\ol P\in\Xi_{k+1}$
        the graph $\Gamma_P$ on the zeros of $P$ 
        with the separatrices as edges. 
    \end{enumerate}
\end{deff}

\begin{prop}\label{prop:graphe sep} 
    The boundary of a connected component of $\Cc\setminus\Gamma_P$
    has either one or two connected components. Each component of the the 
    boundary contains at least one singular point. If the vector field
    is generic, then each component contains at most one singular point.
\end{prop}

This motivates the following definition.

\begin{deff} A connected component of $\Cc\setminus \Gamma_P$ is said to be
    \begin{enumerate} 
        \item a \emph{sepal zone} if its boundary has two components;
        \item a \emph{petal zone} if its boundary has one component.
    \end{enumerate}
\end{deff}

Here is the motivation for this terminology. A petal zone contains a petal of
the parabolic point at infinity and is bounded by a single curve formed by the
union of separatrices and singular points. A sepal zone is a domain isomorphic to
a strip that fits in the interstice spaces between petals. Note that in the
rectifying coordinate $t = \int P$, a sepal zone becomes a horizontal strip, and
petal zone, an upper or lower half-plane.

\begin{figure}[bh!]
    \includegraphics{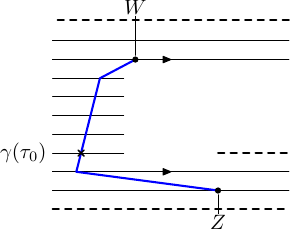}
    \caption{Points $Z,W$ with $\omega$-limit in different ends,
    and three boundaries in dashed lines}
    \label{fig:2 bords}
\end{figure}

\begin{proof}[Proof of Proposition~\ref{prop:graphe sep}]
    Let $X$ be a connected component of $\Cc\setminus\Gamma_P$.
    It is clear that $\del X$ is a union of separatrices and singular points, since
    its boundary is contained in $\Gamma_P$. Moreover, $X$ must be unbounded, since
    $\Gamma_P$ contains no cycle. We call an \emph{end} of $X$ a connected region
    of $X\setminus B(0,R)$ for a very large $R$. We want to prove that $X$ has
    either one end (then $X$ is a petal) or two ends (then $X$ is a sepal).

    Suppose $X$ has more than two ends. First, we transform $X$ in the
    rectifying coordinate $t = \int^z P$, where $\ol P\delz$ becomes
    $|P|^2\delt$. Note that $t$ extends to the boundary into a continuous
    function. Let $Z,W\in t(X)$ be two points such that their $\alpha$-limit or
    their $\omega$-limit goes in different ends. We can of course suppose that
    their $\omega$-limit goes in different ends, otherwise it suffices to
    inverse time. We can join $Z$ and $W$ by a polygonal arc $\gamma\colon [0,1]
    \to \Cc$ transverse to the flow. See Figure~\ref{fig:2 bords}. Let 
    $$
        \tau_0 = \sup \{ \tau\in [0,1]\mid \text{the $\omega$-limit of the orbit of
            $\gamma(\tau)$ is in the same end as $Z$}\}
    $$
    Then the orbit of $\gamma(\tau_0)$ must land on the boundary so that the orbit
    of $t^{-1}(\gamma(\tau_0))$ must land on one the boundaries, a contradiction.
\end{proof}

\subsection{Structurally stable vector fields}
We now specialize our study to the case of structurally stable vector fields,
from here until Section~\ref{sec:non stable}.

Broadly speaking, a vector field is structurally stable if its phase portrait
remains topologically the same after a small perturbation. In the case of
antiholomorphic polynomial vector fields, this is equivalent to having only simple singular points and no
heteroclinic connection. (Limit cycles and homoclinic
loops were already impossible in these systems.)

\begin{notation} We note $\Xi_{k+1}'$ the subspace of \emph{generic} 
    monic centred anti-polynomial vector fields in $\Xi_{k+1}$.
\end{notation}

For structurally stable vector fields, the understanding of the phase portrait can be translated solely in combinatorial terms. The goal of the following discussion is to construct this so called combinatorial data.

\begin{figure}
    \subfigure[Separatrices (incoming in solid, outgoing in dashed) and singular
    points in $\ol \Cc$%
    \label{fig:portrait_de_phase_a}]{\includegraphics{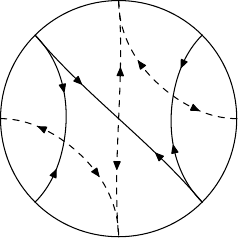}}
    \hskip3cm
    \subfigure[Incoming graph\label{fig:portrait_de_phase_b}]{\includegraphics{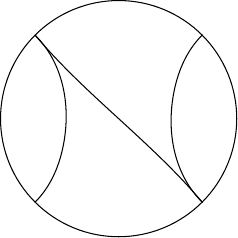}}
    \caption{Incoming graph obtained from a separatrix graph for $k = 2$.}
    \label{fig:portrait_de_phase}
\end{figure}

At each zero $z_j$ of $P$, there are $2$ incoming and $2$ outgoing separatrices (see
Figure~\ref{fig:portrait_de_phase_a}). These come in alternating order since every singularity is a saddle point. 

Given a zero $z_j$ of $P$, define an edge as the union of $z_j$ and its two
incoming separatrices. Now, embedding these edges in $\ol \Cc$ taking the
endpoints as vertices, this gives a graph embedded in $\ol \Cc$, with $k+2 $
vertices and $k+1$ noncrossing edges (see Figure~\ref{fig:portrait_de_phase_b}).
Call this graph the incoming graph of $P$. This data organizes the topological
behaviour of the phase portrait since singular points are saddle points and
incoming separatrices determine the basin of attraction at $\infty$ (see Figure
\ref{fig:portrait_de_phase}). This will be proved in
Theorem~\ref{theo:quasi-conforme}.

We'll refer to this
graph as the \textit{combinatorial invariant} of $\ol P\delz$. It would have
been equivalent to consider edges as union of zeroes and outgoing separatrices instead,
which we call the \emph{outgoing graph}, as the combinatorial invariant. Indeed, 
we have the following result.

\begin{lemma}
    The dual of an incoming graph $G$ is isomorphic to the outgoing graph.
    Moreover, the separatrix graph is given by the union of the incoming graph
    and the outgoing graph in $\ol \Cc$.
\end{lemma}
\begin{proof}
    Recall that $\Cc\setminus G$ is the basin of attraction of infinity.
    If two components share a boundary, then there is singular point
    in this boundary and an outgoing separatrix going in each component,
    forming an edge of the outgoing graph. If $G'$ denote the outgoing graph,
    we can do a symmetric argument on $\Cc\setminus G'$.
\end{proof}

We now want to understand the structure of those graphs. 

\begin{lemma}\label{lem:inc1}
    If $G$ is an incoming graph, then it is embedded (i.e. no two vertices intersect).
\end{lemma}
\begin{proof}
    This follows from Cauchy's uniqueness theorem. Two crossing edges would result from crossing separatrices.
\end{proof}
\begin{lemma}\label{lem:inc2}
    If $G$ is an incoming graph, then $G$ has no cycles.
\end{lemma}

\begin{proof}
    By contradiction, assume that there exists a cycle $c$ and take any edge $e$ in the cycle. The edge $e$ corresponds to the two separatrices entering some zero $z_j$. Let $U\subseteq \ol \Cc$ be the region bounded by the edges of the cycle $c$. Recall that there must exist an outgoing separatrix $\gamma$ at $z_j$ that satisfies $\lim_{t\to\infty}|\gamma(t)| = \infty$. But the only ends of the region $U$ correspond to directions going away from infinity, a contradiction (see Figure \ref{fig:lemme_cycle}).
    \begin{figure}
        \includegraphics{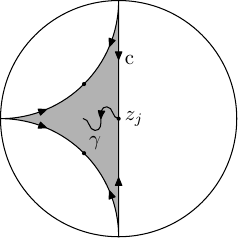}
        \caption{A cycle in an incoming graph prevents any outgoing separatrix at $z_j$ to reach the circle at infinity. } \label{fig:lemme_cycle}
    \end{figure}   
\end{proof}
By combining the previous lemmas, we can count the number of sepal zones
and petal zones in $\Cc\setminus \Gamma$.

\begin{prop}\label{prop:nb de composantes} 
    Let $\ol P\delz$ be a generic vector field of degree $k+1$. 
    The separatrix graph $\Gamma_P$ divides $\Cc$ into $k$ sepal zones
    and $2k+4$ petal zones. 
\end{prop}
\begin{proof}
    An incoming graph $G$ has $k+1$ edges, so it divides $\ol\Cc$ in $k+2$
    connected, simply-connected components. In each component, there is a vertex
    of the dual tree $G^*$ on the circle at infinity. If $n_j$ is the incidence
    of the $j$-th vertex of the dual, it means this connected component gets
    divided into $n_j+1$ connected components by the edges of the dual tree.
    Since the dual is the outgoing tree, it also has $k+1$ edges, so that $\sum_{j=1}^{k+2}
    n_j = 2(k+1)$. Therefore, the number of connected components of
    $\Cc\setminus \Gamma_P$ is $\sum_j (n_j+1) = 2(k+1) + k+2= 3k+4$. 
    
    \vskip3pt
    Moreover, for each $n_j$, the corresponding connected component is divided
    in two petal zones and $n_j-1$ sepal zones. It follows that there are $\sum_j
    (n_j-1) = k$ sepal zones and $2k+4$ petal zones. 
\end{proof}

\section{Combinatorial invariant}\label{sec:combinatoire}

By Lemmas~\ref{lem:inc1} and~\ref{lem:inc2}, the incoming and outgoing graphs
are part of the family of graphs known as noncrossing trees in the literature.

\begin{deff}
    A \textit{noncrossing tree of order $n$} (or $n$-nc tree for short) is a rooted
    tree on $n$ vertices and $n-1$ edges that is embedded in
    $\overline{\mathbb{D}}$ ($\cong\ol\Cc$), where $\mathbb D$ is the unit disk in $\Cc$,
    such that every vertex belongs in
    $\partial\mathbb{D}$ and no two embedded vertices cross.
\end{deff}

We ask that an isomorphism of $n$-noncrossing trees  fixes the root and preserves the cyclic order of the vertices in the boundary circle.

Denote by $\sim_{\text{top}}$ the relation of topological orbital equivalence 
fixing the marked point at infinity (see Definition~\ref{def:eq-top-orb}) on
$\Xi_{k+1}$. Recall that $\Xi_{k+1}'$ is the subspace of monic centred generic
anti-polynomial vector fields.

\begin{theo}\label{theo: bijection}
    For each $k\geq 0$, there is a bijection between $(k+2)$-nc trees and the
    equivalence classes of $\Xi_{k+1}'/{\sim_{\text{\rm top}}}$.
\end{theo}

It will be a corollary of Theorem \ref{theo:realisation} of Section \ref{sec:realisation} that every nc tree can be realized as the combinatorial invariant of some generic antiholomorphic polynomial vector field. Hence, counting the number of orbital topological equivalence classes amounts to counting the number of
noncrossing trees.

\subsection{Counting the number of noncrossing trees}\label{subsec:counting the number of nct}

This section is devoted to counting the number of isomorphism classes of noncrossing trees of order $n$ rooted at $1$. Fix a number of vertices $n\geq 1$. Given a $n$-noncrossing tree, we index every vertex counterclockwise by 1 to $n$ (mod $n$), where $ j$ (mod $n$) represents vertex $e^{(j-1)\frac{2\pi i}{n}}$. We define
\begin{gather*}
    A(n):=\#\{\text{$n$-noncrossing trees}\}.
\end{gather*}

Our main interest is to compute a formula for $A(n)$. This is done more easily by introducing the intermediate quantity 

\begin{gather*}
        A_1(n):=\#\{\text{$n$-noncrossing trees such that $\deg{1} = 1$}\}.
\end{gather*}

and to link $A(n)$ with it.

\begin{figure}
    \includegraphics{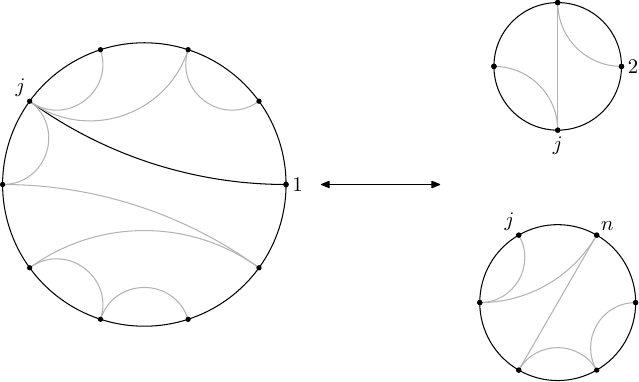}
    \caption{There is a bijection between $n$-nc trees such that $\deg 1 = 1$ having $(1,j)$ as an edge and the set of pairs $(G_1,G_2)$ of nc trees on $j-1$ and $n+1-j$ vertices.} \label{fig:a1}
\end{figure}

\begin{lemma}\label{prop:A_1rec}
    The numbers $A_1(n)$ satisfy $A_1(1) = 0$ and 
    \begin{equation}\label{eq:2}
        A_1(n) = \sum_{j=2}^nA(j-1)A(n+1-j)
    \end{equation}
    for $n\geq 2$.
\end{lemma}
\begin{proof}
    To count $A_1$, remark that every nc tree that satisfies $\deg{1} = 1$ must
    have one and only one edge $e$ coming out of $1$. Let $2\leq j\leq n$ be the
    other end of $e$. Consider the two subgraphs $G_1$ with vertices $2$ to $j$
    and $G_2$ with vertices $j$ to $n$ and set $n_1 := j-1$ and $n_2 := n+1-j$.
    There are no edges from $1\leq m < j$ to $j < \ell \leq n$ and no other edges
    to 1, so $G_i$ must have $n_i-1$ edges, therefore it is a $n_i$-nc tree.
    
    Conversely, suppose we're given two nc trees $G_1$ and
    $G_2$ of order $n_1$ and $n_2$ that satisfy $n_1 + n_2 = n$. 
    Let $j := n_1 +1$. Take the disk
    with $n$ evenly spaced vertices labelled from $1$ to $n$ on it's boundary
    and draw the edge $(1,j)$. We can reconstruct a nc tree of order $n$ by
    mapping the vertices of $G_1$ to the vertices $2$ to $j$ on $S^1$ and the vertices of $G_2$ to vertices $j$ to $n$ on $S^1$ (see
    Figure~\ref{fig:a1}). The resulting graph is a nc tree that satisfies $\deg{1} = 1$.
    
    Since there are $A(n_i)$ such $G_i$, this gives $$A_1(n) = \sum_{j=2}^n A(j-1)A(n+1-j).$$
\end{proof}

\begin{lemma}\label{lem: a} The numbers $A(n)$ satisfy $A(1) = 1$ and
    \begin{equation}\label{eq:3} A(n) = \sum_{j = 2}^n A_1(n+2-j)A(j-1)
    \end{equation} for $n\geq 2$.
\end{lemma}

\begin{deff}\label{def:dualgraph}
    The dual of nc tree $G$ is the usual dual of a planar graph, with the additional
    constraint that the vertices lie at position $e_j' := e^{-2i\pi\over 2k} e_j$ on $S^1$.
    We note it $G^*$.
\end{deff}

For a $n$-nc tree $G$, there are $n$ connected components of $\ol \Dd\setminus G$, hence its dual $G^*$ has $n$ vertices. Every edge of $G^*$ must cross one and only one edge of $G$, therefore $G^*$ has $n-1$ edges, which implies that $G^*$ is a tree. Furthermore, $G^*$ is noncrossing since the dual of a planar graph is planar.

\begin{figure}
    \includegraphics[scale=.75]{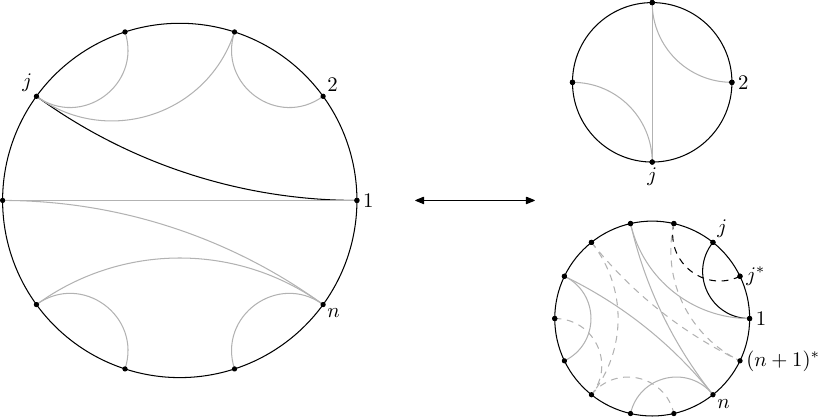}
    \caption{There is a bijection between $n$-nc trees such that $(1,j)$ is an edge and $(1,\ell)$ is not an edge for all $2\leq \ell\leq j-1$ and the set of pairs $(G_1,G_2)$ where $G_1$ is a $(j-1)$-nc tree and $G_2$ is a $(n+2-j)$-nc tree that has $(1,j)$ as an edge. Furthermore, nc trees such as $G_2$ are characterized by having a dual (dashed) with $\deg j^* = 1$.} \label{fig:a}
\end{figure}

\begin{proof}(of Lemma \ref{lem: a}) For each $2\leq
    j\leq n$, we want to count the number of nc trees such that $(1,\ell)$ is
    not an edge for all $2\leq \ell\leq j-1$, but $(1,j)$ is an edge. This
    partitions the set of $n$-nc trees in $n-1$ disjoint subsets (indexed by $j$). Hence, we only
    need to count the number of nc trees in each of those subsets.
    
   Given a vertex $2\leq j\leq n$ such that $(1,j)$ is an edge, if $(1,\ell)$
   is not an edge for all $2\leq\ell\leq j-1$, then we can divide the nc tree into
   $G_1$, the subgraph induced by vertices $2$ to $j$ and $G_2$, the one induced
   by vertices $j$ to $1$, so that $|V(G_1)| = j-1$ and $|V(G_2)| = n+2-j$. Then
   $G_1$ is a $(j-1)$-nc tree and $G_2$ is a $(n+2-j)$-nc tree such that $(1,j)$
   is an edge (see Figure \ref{fig:a}).
     
   The number of such $G_1$ is $A(j-1)$. To count the number of $(n-j+2)$-nc trees
   $G_2$ that have $(1,j)$ as an edge, consider the dual tree $G_2^{*}$. Label
   its vertices by $j^*, (j+1)^*,...,(n+1)^*$ counter clockwise, starting with the
   vertex between $1$ and $j$ (i.e.~$\ell^*$ is the $(\ell-j)$-th vertex at position $\exp\left({
   (\ell-j +1)\frac{2\pi i}{n_2} - \frac{\pi i}{n_2}}\right)$, where $n_2 =
   n+2-j$). Then, to ask for $(1,j)$ to be an edge in $G_2$ is equivalent to
   ask $G_2^*$ to be a $(n+2-j)$-nc tree with $\deg(j^*) = 1$. Hence, there are
   $A_1(n+2-j)$ possibilities for the tree $G_2$. This gives Equation \eqref{eq:3}.
\end{proof}

\begin{theo}
The numbers $A(n)$ of $n$-nc trees satisfy
\begin{equation}\label{eq:4}
    A(n) = \frac{1}{2n-1}\binom{3n-3}{n-1}
\end{equation}
for $n\geq 1$, where $n$ is the number of vertices.

The numbers $A_1(n)$ satisfy $A_1(1) = 0$ and 
    \begin{equation}\label{eq:5}
        A_1(n) = \frac{2}{3n - 4}\binom{3n-4}{n-2}.
    \end{equation}
\end{theo}
\begin{proof}
    For $n\geq 0$, let $B(n) = A(n+1)$. We have
    \begin{equation}\label{eq:B conv}
    \begin{aligned}
        B(n)& = \sum_{j = 0}^{n-1}A_1(n+1-j)A(j+1)\\
        &= \sum_{j = 0}^{n-1}\sum_{k=0}^{n-j-1}A(k+1)A(n-j-k)A(j+1)\\\noalign{\penalty0}
        &=\sum_{j = 0}^{n-1}\sum_{k=0}^{n-j-1}B(k)B(n-j-k-1)B(j)\\
        &=\mathop{\sum_{a+b+c = n-1}}_{0\leq a,b,c\leq n-1} B(a)B(b)B(c).
    \end{aligned}
    \end{equation}
    Then if $\mathcal{B}(z)$ is the generating function, the recurrence relation
    gives $z\mathcal{B}^3(z) = \mathcal{B}(z) - 1$. A consequence of Lagrange's Inversion
    Formula states that a formal series $\mathcal A(z) = \sum_{n=0}^\infty a_n z^n$
    satisfying $\mathcal A(z) = a_0 + z g\bigl( \mathcal A(z)\bigr)$ will have coefficients
    $$
        a_n = {1\over n!} {\d^{n-1} \over \d z^{n-1}}\!\! \left.\Bigl( g(z)^n\Bigr)\right|_{z=a_0}
        \rlap{\qquad for $n\geq 1$.}
    $$
    Applying this formula with $a_0 = 1$ and $g(z) = z^3$ to $\mathcal B(z)$ 
    then gives us that $B(n) = \frac{1}{2n+1}\binom{3n}{n}$ from which
    the proposition follows. This is done in \cite{NOY98}.

    The formula for $A_1(n)$ is done in a similar way; see \cite{NOY98} for more details. The sequence $\bigl(
    A_1(n+1)\bigr)_{n\geq 1}$ is \href{https://oeis.org/A006013}{A006013} in the
    OEIS.
\end{proof}

The sequence $A(n)$ is well known and corresponds
\href{https://oeis.org/A001764}{A001764} in the OEIS up to a reindexing. Those
numbers also count ternary trees, a fact that we will use in
Section~\ref{sec:arbre ternaire}. Indeed, let $T(n)$ be the number
of ternary tree with $n$ internal vertices, since there are three subtrees
attached to the root, it follows that $T(n)$ verifies Equation~\ref{eq:B conv},
and since $T(0) = 1$, $T(1) = 1$, we obtain $T(n) = A(n+1)$.

We will see in Section~\ref{sec:realisation} (see Remark~\ref{rem:nb strates})
that $A(k+2)$ counts the number of generic strata in the parameter space for
generic anti-polynomial vector fields of degree $k+1$. By \emph{generic
stratum}, we mean an open connected set in the parameter space for which each
vector field is orbitally topologically equivalent to one another.

Let $A^r(n)$ be the number of $n$-nc trees modulo rotations. Looking at some examples shows us that some nc-crossing tree may be obtained from others by rotation. The goal of the following discussion is to compute the number of nc trees modulo rotation.
\begin{lemma}\label{lem:Rotation}
    Let $G$ be a $n$-nc tree, then $G$ has a non-trivial rotational symmetry $\alpha$ only if $n$ is even. Moreover, if this is the case, $\alpha$ necessarily fixes a unique edge with antipodal endpoints.
\end{lemma}
\begin{proof}
    Suppose $\alpha$ acts on $G$ by rotation (i.e. it translates every vertex by a constant amount mod $n$). Index every vertex by $0$ to $n-1$ (mod $n$). Let $m = \ord\alpha$. Then we have $m|n$.
    We first claim that there exists some subgraph $H$ of $G$ such that $G = \bigcup_{0\leq j\leq m-1}\alpha^jH$. 
    Indeed, let $d = \tfrac n m$ and let $S = \{0,...,d-1\}$. Define $H$ as the subgraph induced by the edges that have an endpoint in $S$. Then $H$ satisfies the claim.
    
    If $f\in E(G)$ is an edge, we say that $f$ is \textit{internal} if both its endpoints are in $S$ and \textit{external} if it has only one endpoint in $S$. Denote by $e$ the number of external edges of $H$ and by $i$ the number of internal edges of $H$. Counting the edges over each $\alpha^j H$, internal edges are counted only once, but external edges end up in two different graphs $\alpha^{j_1}H$ and $\alpha^{j_2}H$. Hence, the total number of edges in $G$ is 
    \begin{equation}\label{eq:arete}
        n-1 = mi+m\frac{e}{2}.
    \end{equation}
    This gives that $m|2(n-1)$. Hence we must have that $m|2$, but since $\alpha$ is not trivial, we must have that $m=2$, thus $n$ is even and $\alpha$ has order 2.
    
    Now, suppose that $\alpha: j \mapsto j+\frac{n}{2}$ (mod $n$) is a rotational symmetry of $G$ and let $H$ be the subgraph as above. Remark that there must exist an external edge $f$ since, otherwise, we would have that $\#\{\text{edges of $G$}\} = n-1 = 2i$, contradicting the parity of $n$.

    Let $\pi:\mathbb{Z}\to \mathbb{Z}/n\mathbb{Z}$ be the projection and lift $G$ to $\mathbb{Z}$. Conjugating $G$ by a rotation, we may assume that $\tilde f
    := (0,j)$ is a lift of $f$, with $\frac{n}{2}\leq j\leq n-1$. The rotation $\alpha$ lifts to a translation $\tilde \alpha$ by $\frac n 2$, so that $\tilde \alpha.\tilde{f}= (\frac{n}{2},j+\frac{n}{2})$. 
    Now, whenever $\frac n 2 < j\leq n-1$, the edge $\alpha.f = \pi(\tilde \alpha. \tilde f)$ and $f$ cross in $G$, so we must have that $j=\frac n 2$, meaning that $f$ has antipodal endpoints. Since the graph is noncrossing, this edge must be unique.
\end{proof}

\begin{theo}
    Denote by $A^r(n)$ the number of $n$-nc trees modulo rotation. We have that
    \begin{equation}\label{eq:6}
        A^r(n) = \frac{1}{n}\bigg(A(n) + \mathbf{1}_{n\text{ even}} {n \over 2} A_1\bigg(\frac{n}{2} + 1\bigg)\bigg)
    \end{equation}
\end{theo}
This is sequence \href{https://oeis.org/A296532}{A296532} in the OEIS up to a reindexing. This result
was probably already proved, but we could not find a reference, therefore, we
added a proof here. A discussion about similar cases may be found in \cite{NOY98}.
\begin{proof}
    Let $G$ be $n$-nc tree with vertices $0\leq j\leq n-1$.
    
    In the case where $n$ is odd, Lemma \ref{lem:Rotation} tells us that there are precisely $n$ elements in every rotation class,
    so that $nA^r(n) = A(n)$. 

    When $n$ is even, every nc tree that has a rotational symmetry necessarily has $m := \frac{n}{2}$ elements in its rotation class. Therefore, to compute $A^r (n)$, we need to count the number of nc trees that have a rotational symmetry to artificially add them up before dividing by $n$. 
    To count such $G$'s, by Lemma \ref{lem:Rotation}, there exists $\alpha$ be the rotation of order $2$, and $f = (0,m)$ its fixed edge. Since $f$ splits $G$ in two $(m+1)$-nc trees $G_1$ and $G_2$ and that those two trees mutually determine the other one via $\alpha$, we're left to compute the number of $(m+1)$-nc trees $G_1$ that have $(0,m)$ as an edge. Using the same argument as in Lemma \ref{lem: a} and considering the dual tree tells us that this number is $A_1(m+1)$. Therefore, we obtain
    $$
        A(n) = nA^r(n)  - mA_1(m + 1)
        \rlap{\qquad (for $n$ even).}
    $$
\end{proof}

\subsection{Characterization of the topological phase portrait with noncrossing
    trees}

\begin{deff}
    Given a generic polynomial $\ol P \in\Xi'_{k+1}$, the \emph{combinatorial invariant} 
    of $\ol P\delz$ is the $(k+2)$-nc tree induced by the separatrix graph, with root at $\exp\bigl({ {2 \pi i\over 2(k+2)}}\bigr)$ on the circle at infinity.
\end{deff}

\begin{theo} \label{theo:quasi-conforme}
    Let $\ol{P_1},\ol{P_2} \in \Xi'_{k+1}$. If $\ol{P_1}\delz$ and $\ol{P_2}\delz$
    have the same combinatorial invariant, then they are quasi-conformally equivalent.
    (In other words, there exists a quasi-conformal mapping $h$ that maps orbits of
    $\ol{P_1}\delz$ on orbits of $\ol{P_2}\delz$). In particular, the two vector fields are 
    topologically orbitally equivalent.
\end{theo}

\begin{proof}
Let $C_i=\Cc\setminus \Gamma_{P_i}$ where $\Gamma_{P_i}$ is the separatrix graph. Since $P_1$ and $P_2$ have the same
combinatorial invariant, there is a natural correspondence between
the connected components $X_1^j$ of $C_1$ and $X_2^j$ of $C_2$. By Proposition \ref{prop:graphe sep}, each $X_i^j$ is
adjacent to either one or two singular points. In the first case, we denote it by $z_i^j$ and in the second case, by $z_i^j$ and $\tilde z_i^j$. When there are two singular points, we
label them so that $\Im \int_{z_i^j}^{\tilde z_i^j}
P(\zeta) d\zeta > 0$. (When $P$ is generic, the latter is always non zero; we will discuss this in Remark \ref{rem:orientation segment}.)
    
Define $t_i^j(z) = \int_{z_i^j}^zP_i(\zeta) d\zeta$. Then the pushforward of
$\ol{P_i}\delz$ by $t_i^j$ gives 
\begin{align*} \ts
    {t_i^j}_*\bigl(\overline{P_i}\frac{\partial}{\partial z}\bigr) = \ds
    |P_i|^2\frac{\partial}{\partial t_i^j},
\end{align*}
    a real scaling of the horizontal vector field $
    1\frac{\partial}{\partial t_i^j}$.
    
    If $X_i^j$ is a sepal zone, this
    gives a conformal equivalence between $X_i^j$ and an infinite open strip
    which we denote $B_i^j$. If $X_i^j$ is a petal zone, then $t_i$ maps it
    conformally to a half-plane that we denote by $H_i^j$. In any case, $t_i^j$
    extends to the adjacent separatrices, giving a homeomorphism from
    $\overline{X_i^j}$ onto $\overline{B_i^j}$ or $\overline{H_i^j}$.
    
    If $X_1^j$ is a petal zone that corresponds to $X_2^j$,
    define $\psi_j:H_1^j\to H_2^j$ as identity map. If $X_1^j$ and
    $X_2^j$ are sepal zones, let $\tau_i^j = t_i^j(\tilde z_i^j) -
    t_i^j(z_i^j)$, so that $B_i^j = \{a+b\tau_i^j:a\in\mathbb{R},\,b\in(0,1)\}$. We define
    $$\begin{aligned}
        \psi_j:B_1^j &\to B_2^j\\
        a+b\tau_1^j &\mapsto a + b\tau_2^j.
    \end{aligned}$$

    Recall that the rectifying coordinate $t$ maps orbits to horizontal lines. Since the $\psi_j$'s we've just built map horizontal lines to horizontal lines, the composition $h_j :=(t_2^j)^{-1}\circ\psi_j\circ t_1^j$ maps orbits $\ol P_1\delz$ in $X_1^j$ to orbits of $\ol P_2\delz$ in $X_2^j$.
    
    Denote by $\mu_j(z) = \frac{\overline{\partial}\psi_j(z)}{\partial \psi_j(z)}$ the Beltrami coefficient of $h_j$. One can compute that $\mu_j(z) = \frac{\tau_2-\tau_1}{\overline{\tau_1}- \tau_2}$ and the condition $|\mu_j(z)|<1$ is seen to be equivalent to $4\Im(\tau_1)\Im(\tau_2)>0$ which is true. Note that all $\psi_j$'s are invertible and extend to homeomorphisms between $\overline{A_1^j}$ and $\overline{A_2^j}$, where $A_i^j$ is either a strip or a half-plane. Hence, we have the following diagram for every $j$.

\[
\begin{tikzcd}
        X_1^j \arrow{r}{h_j} \arrow{d}{t_1^j} & X_2^j \arrow{d}{t_2^j} \\
        A_1^j \arrow{r}{\psi_j}[swap]{q.c.} & A_2^j         
\end{tikzcd}
\]

Since $t_1^j$, $\psi_j$ and $t_2^j$ extend to the separatrices, so does $h_j$, so we can glue the $h_j$'s in a global homeomorphism $\Psi:\Cc\to\Cc$ that maps orbits of $\ol{P_1}\delz$ to orbits of $\ol{P_2} \delz$. We know that $\Psi$ is quasi-conformal on $\Cc\setminus \Gamma_{P_1}$, hence is quasi-conformal on the whole plane $\Cc$.
\end{proof}

Let $\rho(z) = \lambda z$ be a rotation of order $k+2$. It acts on
$(k+2)$-noncrossing trees by rotating them when seen as embedded in
$\ol{\mathbb{D}}$. A direct consequence of the the preceding theorem is the following.
\begin{cor}
   Suppose that $\overline P_j{\delz}$ has combinatorial invariant
   $G_j$, for $j = 1,2$. If there exists a rotation $\rho$ such that $\rho.G_1 =
   G_2$, then $\rho_*\left(\overline P_1 \delz\right)$ is topologically
   orbitally equivalent to $\overline P_2 \delz$.
\end{cor}

As a consequence, there is $A^r(k+2)$ generic strata in $\Xi_{k+1}'$
modulo rotation.

\def\line(#1 #2 #3){(#1 \ \ #2 \ \ #3)}%
\def\V{v}%
\section{Analytic invariant}\label{sec:analytique}
The analytic invariants of both ${1\over P}\delz$ and $\ol P\delz$ will be given
in terms of $\int P$. Indeed, this corresponds to both the complex time of
${1\over P}\delz$ and the circulation and the flux of $\ol P\delz$.
We will define below the invariant as integrals across the sepal zones,
so we will establish an order on the zones that depends only on the combinatorial
invariant. 

\begin{figure}[h]
    \includegraphics[scale=.8]{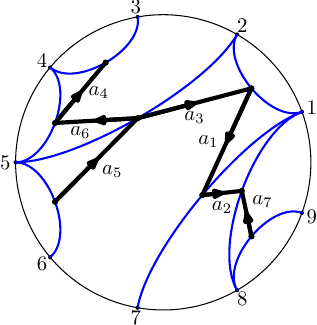}
    \caption{Line segments $a_j$ for the analytic invariants. Notice that at the vertex 5,
    $a_5$ is on the right and $a_6$, on the left, because of the partial order on the edges.}
    \label{fig:segment}
\end{figure}

\begin{deff}
    \label{deff:ordre}
    Let $\ol P\delz \in \Xi_{k+1}'$ and let $G$ be its $(k+2)$-nc tree. Recall that the vertex~$j$ (mod $k+2$) is at position $\exp\!\left( {2i\pi(2j-1)
    \over 2(k+2)} \right)$ on the circle.
    \begin{enumerate}
        \item We define a partial order on the edges. We say $\{j,m\} < \{j,\ell\}$
            if $j < m < \ell$, with $m,l\in\{j+1,j+2,...,j+k+1\}$. In other words,
            the edges attached to $j$ are ordered starting from $j$ and going around the circle
            anti-clockwise.
        
        \item We define an order on the sepal zones. Let $n_j$ be the incidence of the $j$-th vertex of $G$. Then there
            are $n_j-1$ sepal zones associated to this vertex (see Proposition~\ref{prop:nb de composantes}). For $1\leq m\leq k$, there exist unique $1\leq j\leq k+2$ and $1\leq
            \ell \leq n_j-1$ such that $m = (n_1-1) + \cdots +
            (n_{j-1}-1)+\ell$. This means that the $m$-th sepal zone is the $\ell$-th sepal zone of the $j$-th vertex.
            See Figure~\ref{fig:segment}.
        \item 
            Let $e_1 < \cdots < e_{n_j}$ be the edges attached to $j$ and $z_{e_\ell}$ be
            the singular point on $e_\ell$. Then we define 
            $$
                {\rm segment}(j,\ell) = [z_{e_\ell}, z_{e_{\ell+1}}]
                \rlap{\qquad for $1\leq \ell \leq n_j-1$.}
            $$
            For the $m$-th sepal zone, we define the line segment $ a_m := \text{segment}(j,l)$, where $j$ and $\ell$ are as in (2).
    \end{enumerate}
\end{deff}

The way we choose to order the sepal zones is not intrinsically important for our
study, however, we must choose one to go forward.

\begin{deff}[Analytic invariant]\label{def:invariant analytique} Let $\ol P\in
    \Xi_{k+1}'$, $G$ be the combinatorial invariant and $a_j$ be the line segments
    defined in Definition~\ref{deff:ordre} using $G$. The \emph{analytic invariant}
    of $\ol P\delz$ and ${1\over P}\delz$ is the vector $\eta\in \Hh^k$, with the $j$-th
    component $\eta_j = \int_{a_j} P(z)\d z$. 
\end{deff}

\begin{rem}\label{rem:orientation segment} 
    In the definition, we assumed that $\Im \eta_j > 0$. This is indeed the case
    because of the choice of orientation of the $a_j$'s. We may suppose that
    $a_j$ is contained in its sepal zone and make a constant angle in $(0,\pi)$
    with the vector field. In other words, the canonical normal vector of $a_j$
    in the rectifying coordinate will make a constant angle in $(-{\pi\over
    2},{\pi\over 2})$ with the flow, so that the flux of $\ol P\delz$ across
    $a_j$ will be positive, and this is precisely the imaginary part of
    $\int_{a_j} P(z)\d z$. 
%
\end{rem}


The following result proves that the analytic invariant is indeed
invariant and, together with the combinatorial invariant, it describes
the whole dynamics.
\begin{theo} Let $k\in \N$ and $P_1,P_2$ be two monic centred polynomials of
    degree $k+1$. If the two vector fields $\ol P_1\delz$ and $\ol P_2\delz$ are
    generic and have the same combinatorial and analytic invariants, then they are
    analytically equivalent. 
\end{theo}
\begin{proof}
    We obtain a mapping $h\colon \Cc\to \Cc$ such that $h_*\bigl({1\over
    P_1}\delz\bigr) = {1\over P_2}\delz$ using the same construction as in the
    proof of Theorem~\ref{theo:quasi-conforme}, but now the mappings between the
    horizontal strips in the rectifying coordinates are replaced with the
    identity, since the analytic invariants are the same (i.e.~the corresponding
    strips are equal). The rectifying coordinates and the identity are
    holomorphic and extend to the separatrices continuously, which defines~$h$
    as holomorphic bijection.
%
\end{proof}

Recall that a rotation $\rho(z) = \exp\!\left(2i\pi j\over k+2\right)$ acts
on $(k+2)$-nc trees by a rotation of the vertices $\ell \mapsto \ell+j$ (mod $k+2$)
and of the edges. Such a rotation will also change the order of the components of the analytic invariant
by translating the indices $m \mapsto m+(n_{k+2-j+1}-1)+\cdots+({n_{k+2}}-1)$
(mod $k$).

\begin{cor} Suppose a generic $\ol{P_j} \delz$ has invariants $(G_j,\eta^{(j)})$
    for $j=1,2$.
    If there exists a rotation $\rho$ that acts on $(G_1,\eta^{(1)})$
    so that $\rho.G_1 = G_2$ and $\rho_{G_1}.\eta^{(1)}
    = \eta^{(2)}$, then $\rho_*\bigl(\ol{P_1} \delz\bigr)$ and
    $\ol{P_2} \delz$ are analytically equivalent, where, with $\rho(z) = \lambda z$,
    $$
        \rho_*\bigl(\ol P_1\delz\bigr) = \lambda \ol{\ts P_1\bigl({z\over 
        \lambda}\bigr)}.
    $$
\end{cor}
\begin{proof} 
    If a rotation acts on $G_1$, then it must be a rotation of order $k+2$,
    therefore $\rho_*\bigl(\ol{P_1}\delz\bigr)$ is a monic centred polynomial vector field. Its
    invariants are precisely $\rho.G_1$ and $\rho_{G_1}.\eta^{(1)}$, so that it 
    must be equal to $\ol{P_2}\delz$.
\end{proof}

\section{Realization}\label{sec:realisation}
This is the main result of this section.

\begin{theo}[Realization] \label{theo:realisation}
    Let $k\in\N$. For every pair $(G,\eta)$, where
    $G$ is a $(k+2)$-noncrossing tree and $\eta\in\Hh^k$, there
    exists a monic centred generic polynomial $P$ of degree $k+1$ such that $\ol P\delz$
    has $(G,\eta)$ for its invariant. 
\end{theo}

\sloppy
\begin{rem}\label{rem:nb strates} Let $\joli A_k = \{\text{($k+2$)-noncrossing trees}\}$. The moduli
    space of the generic vector fields of $\Xi_{k+1}$ is then $\joli A_k\times
    \Hh^k$, an open set with $\# \joli A_k = A(k+2) = {1\over 2k+3} {3k+3 \choose
    k+1}$ connected components (see Equation~\eqref{eq:4}). We call those components the \emph{generic
    strata}. Each generic stratum of the parameter space can be analytically
    parametrized by $\eta\in\Hh^k$.
\end{rem}
\par\fussy

If we allow homeomorphisms to rotate points at infinity, this corresponds to orientation-preserving 
topological orbital equivalence. In this case, we can quotient $\joli A_k\times \Hh^k$ by the relation $\sim$ defined by
$(G_1,\eta^{(1)}) \sim (G_2,\eta^{(2)})$ if and only if there exists a rotation
$\rho$ of order $k+2$ such that $(G_2,\eta^{(2)}) = (\rho.G_1,
\rho_{G_1}.\eta^{(1)})$.

With a rotation of order $k+2$, we can normalize a vector field $\ol P\delz \in
\Xi_{k+1}$ by asking that $\arg(\ol{P(0)}) \in \bigl[0, {2\pi\over k+2}\bigr)$.
We can choose the invariants of this vector field as a representative of the
class of $[G,\eta]$.

\begin{cor} 
    Let $\widetilde \Xi_{k+1}'$ be the set of generic vector fields $\ol P\delz$
    such that $P$ has degree $k+1$, is monic centred and normalized with
    $\arg(\ol{P(0)}) \in \bigl[0, {2\pi\over k+2}\bigr)$. Then the moduli space
    of those generic vector fields is $\joli A_k\times \Hh^k /{\sim}$. It has
    $A^r(k+2)$ connected component, where $A^r(k+2)$ is given by
    Equation~\eqref{eq:6}.
\end{cor}

The proof of the theorem is similar to the realization construction found
in~\cite{DES} for a generic vector field $P\delz$. The combinatorial data and
the analytic invariant allow us to construct a Riemann surface $\widetilde M$ with open
sets of $\Cc$ in such a way that the vector field $1\delz$ is well-defined on
$\widetilde M\setminus Z$, where $Z$ is a finite set. Then it is proved that $\widetilde M$ has the
conformal type of the Riemann sphere and the push-forward of $1\delz$
will be ${1\over P}\delz$ for some polynomial $P$, after perhaps
applying a Möbius transformation. To obtain $\ol P\delz$ is a bit more
tricky, since the push-forward of an antiholomorphic vector field is not, in
general, antiholomorphic. The trick is to push-forward a well-chosen
log-harmonic function: if $h\colon \widetilde M\to \hat\Cc$ is the isomorphism
with $h_*\bigl(1 {\del\over\del\zeta}\bigr) = {1\over P}\delz$,
then we will have $h_*\bigl(|h'|^{-2}{\del\over\del\zeta}\bigr) = \ol P\delz$.

\subsection{Construction of the Riemann surface}
\label{sec:constr SR}

Consider a $(k+2)$-nc tree $G$ and a vector $\eta\in \Hh^k$. Recall
that the dual $G^*$ is also a $(k+2)$-nc tree. If we embed both $G$
and $G^*$ in the same closed disk $\ol{\mathbb D}$ 
so that the $j$-th vertex of $G$ is mapped to $e_{2j-1}$ and the
$j$-th vertex of $G^*$, on $e_{2j-2}$, where $e_\ell := \exp\!\left(
2i\pi {\ell\over 2(k+2)}\right)$, they divide the
open disk in $3k+4$ connected and simply connected components (see Proposition~\ref{prop:nb de composantes}),
provided the embedding is chosen so that two edges intersect at most
once. We will call the intersection of two edges a \emph{singular
point}, which divides those edges in four arcs that we will call the
\emph{separatrices}. By Proposition~\ref{prop:nb de composantes}, there
are $2k+4$ components adjacent to only one singular point, which we call
the \emph{petal zones}, and $k$ components adjacent to two singular points,
the \emph{sepal zones}.

A petal zone in $\mathbb D$ can be mapped by a homeomorphism to
$\Hh$ or $-\Hh$, which we can choose so that it extends
to the boundary as to map the singular point to $0$; of course,
the separatrices must map on $\R$ on each side of $0$.

The sepal zones are ordered by $G$ from 1 to $k$ according to Definition \ref{deff:ordre}. The $j$-th sepal
zone can be mapped by a homeomorphism to a horizontal strip $B_j$ of width
$\Im\eta_j$, which we can choose in a way that it extends to the boundary
so as to map the singular points, say $z_0$ and $z_1$, on $0$ and $\eta_j$.  There
are two orientations for this; to choose the correct one, note that a sepal
zone has an end on a vertex $e_{\ell^-}$ of $G$ and the other end on a vertex $e_{\ell^+}$ of $G^*$.
We choose a homeomorphism that maps a horizontal line in the strip oriented from left to right
to a curve in $\Dd$ from $e_{\ell^-}$ to $e_{\ell^+}$.
Then the four separatrices in $\Dd$ are forced to be mapped on horizontal lines on each side
of $0$ and $\eta_j$ in the strip.

We start constructing $M$ by sewing the half-planes and horizontal strips along
a part of their boundary if and only if the image of this boundary in $\Dd$ is
the separatrix in $\mathbb D$. At this point, $M$ has the topological type of
$\mathbb D\setminus\{\text{singular points}\}$. The important point is that
the vector field $1\delZ$ is well-defined on $M$. The following discussion
explains how to plug the holes at singular points and define a chart at
$\infty$.

\begin{figure}
    \includegraphics{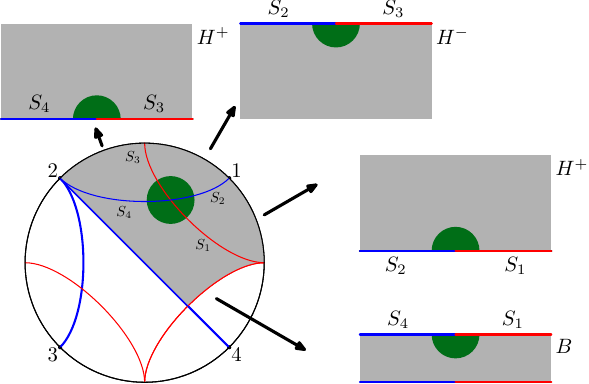}
    \caption{Construction of a holomorphic chart around a singular point. The
    half-disks form a topological disk around a singular point; the $\psi_m$'s
    define a holomorphic chart.}
    \label{fig:carte hol zero}
\end{figure}

\paragraph{Neighbourhood of a singular point}
Let $Z_0$ be a singular point and $S_1,S_2,S_3,S_4$ the four separatrices
attached to it, with $S_1$ an outgoing separatrix (any one of the
two) and $S_2,S_3,S_4$ going anti-clockwise around $Z_0$. Then
$Z_0$ appears on the upper boundary of two charts, and lower boundary of two
charts. For $m\in \{1,2,3,4\}$, we define $D_m$ as a upper (for $m$ odd) or
lower (for $m$ even) half-disk from $S_m$ to $S_{m+1}$ (with $S_5 := S_1$) with
small enough radius centred at $Z_0$ and contained in the chart. If we glue the
$D_m$'s along the separatrices, we obtain a neighbourhood $D_{Z_0}$ of $Z_0$.
Let
$$
    \psi_m(Z) = (Z - Z_0)^{1/2} \quad\text{with $\arg (Z-Z_0) \in [(m-1)\pi,m\pi]$.}
$$
This defines an isomorphism $\psi_{Z_0}\colon D_{Z_0}\setminus \{Z_0\} \to 
D(0,r)\setminus\{0\}$ for some $r>0$. By Riemann's Removable Singularity Theorem, we can extend $\psi_{Z_0}$ to $Z_0$ and obtain
a chart $(D_{Z_0},\psi_{Z_0})$ of $Z_0$. See Figure~\ref{fig:carte hol zero}.
We can verify that $(\psi_{Z_0})_*(1\delZ)(z) = {1\over 2z} \delz$, so that the
dynamics is that of a simple saddle node.

Doing this at every singular point extends $M$ to a Riemann surface homeomorphic to $\mathbb{D}$.

\paragraph{Neighbourdhood of infinity}

To determine a conformal equivalence of a neighbourhood of infinity
on $M$ to a neighbourhood in $\Cc$ is more delicate than the
previous paragraph. Indeed, the expected dynamics is one of a
parabolic point and there are no canonical separatrices to divide
the neighbourhood for us in $\Cc$. A neighbourhood of infinity will
be defined to be any open set of $M$ that contains an annulus $\{r <
z < 1\}$ when seen in $\mathbb D$. See Figure~\ref{fig:collage
infini a}. On the topological level, this correspond to quotienting
$\del\mathbb D$ to a single point.  We now want to define a
holomorphic chart on such a neighbourhood of infinity.

\begin{figure}
    \centering
    \subfigure[Neighbourhood of infinity represented in $\ol\Delta$,
        given a noncrossing tree, with $\joli E_j^\pm$ representing
        the ends of sepal zones\label{fig:collage infini a}]{\qquad\includegraphics{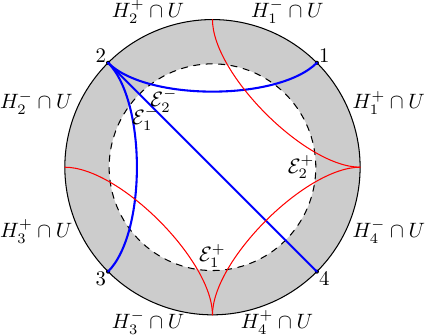}}
    \subfigure[Neighbourhood of infinity in the half-planes glued with
        ends and their translation; the separatrices are named
        $S_0,S_1,\ldots$ to simplify%
        \label{fig:collage infini b}]
        {\centerline{\includegraphics{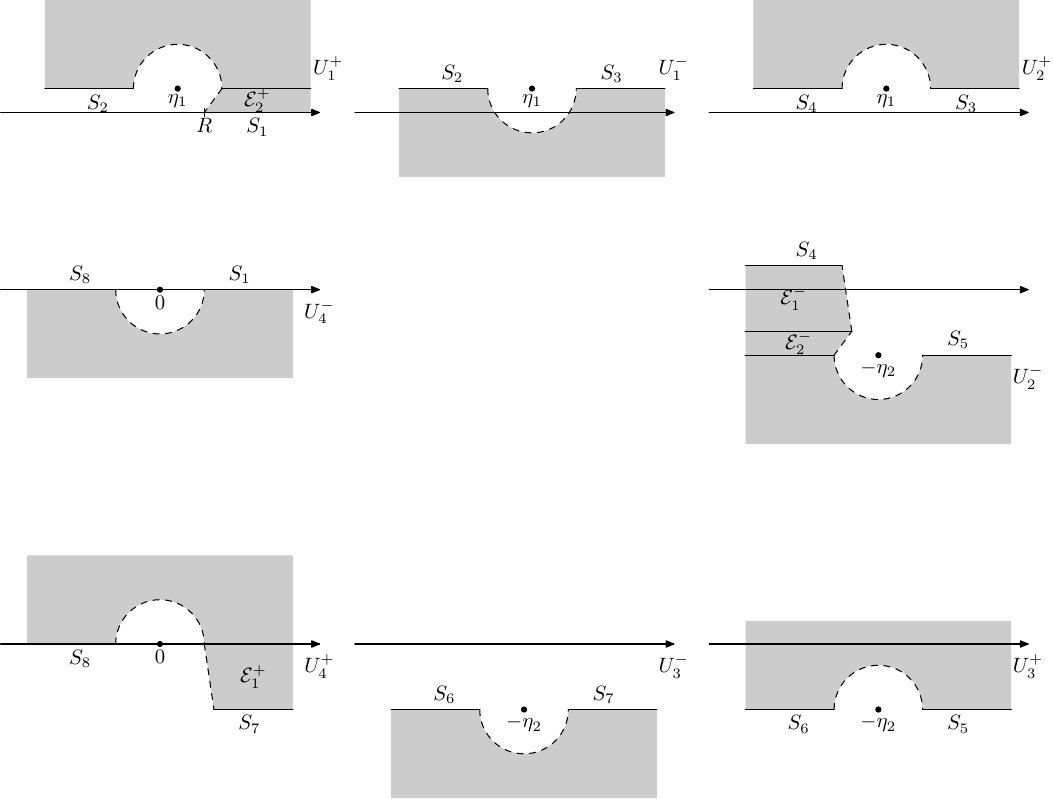}}}
    \caption{Example with $k=2$ of a neighhood of infinity given a 2-noncrossing
        tree and an analytic invariant $\eta = (\eta_1,\eta_2)$}
    \label{fig:collage infini}
\end{figure}

We will label vertices of $G^*$ by $1^+,\ldots, (k+2)^+$, with the $j$-th vertex
positioned at $\exp\!\left( {2i\pi (j-1)\over k+2}\right)$.  Let $U$ be a
neighbourhood of infinity. For a strip $B$, $B\cap U$ has two components that we
call ends. (See Figure~\ref{fig:collage infini a}.) Let $E_j^\pm = \{1\leq m\leq
k \mid B_m\text{ has an end at vertex } j^\pm\}$ (it may be empty).  For each
$m\in E_j^\pm$, we define $\joli E_m^\pm$ as the end of $B_m$ at vertex $j^\pm$,
where the strips $B_1,\ldots,B_k$ are ordered following the sepal zones.

Starting with the petal zone $H_1^+$ between the vertices $1^+$ and $1$, we glue
the positive $\joli E_m^+$ ($m\in E_1^+$) ends at $1^+$ (if any) to $H_1^+\cap U$ and we translate by
$\xi_1^+ = \sum_{m\in E_1^+} \eta_m$, see Figure~\ref{fig:collage infini b}. We obtain an open set $U_1^+$ of $\Cc$. We
define $\psi_1^+ \colon U_1^+ \to \Cc$ by $\psi_1^+(Z) = Z^{-{1\over k+2}}$,
with $\arg(Z) \in \bigl(-{\pi\over 2}, {3\pi\over 2}\bigr)$.

For the next petal zone $H_1^-$ (between $1$ and $2^+$), we glue
negative ends $\joli E_m^-$ ($m\in E_1^-$) at $1$ (if any) to $H_1^- \cap U$ and we translate by 
$\xi_1^- = \sum_{m\in E_1^+} \eta_m - \sum_{m\in E_1^-}\eta_m$ to obtain an
open set $U_1^-$.  We then define $\psi_1^-\colon U_1^- \to \Cc$ by $\psi_1^-(Z)
= Z^{-{1\over k+2}}$ with $\arg(Z) \in \bigl( {\pi\over 2}, {5\pi\over
2}\bigr)$. The translation by $\xi_1^-$ is necessary so that when
$U_1^+$ and $U_1^-$ are glued, the $\psi_1^\pm$ extend to an injective map on
the quotient into a connected open set of $\Cc$.

For the general case, we glue the positive (resp.~negative)
ends at $j$ (resp.~$j^+$) to $H_j^+\cap U$ (resp.~$H_j^-\cap U$)
and we translate by 
$$
    \xi_j^+ = \sum_{\ell=1}^{j-1} \Bigg[ \sum_{m\in E_\ell^+} \eta_m
      - \sum_{m\in E_\ell^-} \eta_m \Bigg] + \sum_{m\in E_j^+} \eta_m
    \quad
    \left( \hbox{resp.~} \xi_j^- = \sum_{\ell=1}^{j} \Bigg[\sum_{m\in E_\ell^+} \eta_m
      - \sum_{m\in E_\ell^-} \eta_m \Bigg] \right).
$$
We obtain an open set $U_j^+$ (resp.~$U_j^-$) of $\Cc$. We define $\psi_{j,\infty}^+
\colon U_j^+\to \Cc$ (resp.~$\psi_{j,\infty}^-\colon U_j^-\to\Cc$)
by $\psi_j^\pm(Z) = Z^{-{1\over k+2}}$ with $\arg(Z) \in 
\bigl((4j-5){\pi\over 2}, (4j-1){\pi\over 2}\bigr)$ (resp.~in
$\bigl((4j-3){\pi\over 2}, (4j+1){\pi\over 2}\bigr)$).
See Figure~\ref{fig:collage infini b}. This defines a conformal map
$\psi_\infty$ on $U$ that sends $U$ on some punctured neighboorhood
of $0$.  We can then add $\infty$ to $M$ and $\psi_\infty$ extends
as a holomorphic chart of $\infty$ on $M$.

\def\Mtilde{\widetilde M}%
\paragraph{Equivalence with the Riemann sphere}
Let $\Mtilde$ be the Riemann surface obtained from $M$ after adding the
singular points $\{Z_j\}_{1\leq j\leq k+1}$ and $\infty$. It is a compact manifold with no
boundary. We can compute the Euler characteristic with the following
triangulation: the vertices are the $k+1$ singular points and
$\infty$, the edges are the $4(k+1)$ separatrices, and the faces are the $2(k+2)$
half-planes and the $k$ strips. We have $k+2$ vertices, $4(k+1)$
edges and $3k+4$ faces, so that $\chi_{\Mtilde} = 2$. Therefore, $\Mtilde$ is a
topological sphere and thus, by Riemann's Uniformization Theorem, it has the
conformal type of the Riemann sphere $\hat \Cc = \Cc\cup \{\infty\}$.

\subsection{Proof of the Realization Theorem}\label{sec:dem realisation}%
We now prove Theorem~\ref{theo:realisation}.

Let $h\colon \Mtilde \to \hat\Cc$ be a conformal mapping. Since we can post-compose 
by a Möbius transformation, we can suppose that $h(\infty) = \infty$. We can further
normalize $h$ by requiring that $h$ in the chart $(\infty,\psi_\infty)$ has the
form $\zeta + o(\zeta)$. Note that $h$ is still not unique, as we can post-compose by an
affine transformation $z\mapsto \lambda z + b$, where $\lambda$ is a $(k+2)$-root of unity.
To fix the rotation, we require that the real interval $(R,\infty)$, for $R$ large enough,
in $H_1$ be mapped by $h$ in the attracting sector at infinity that intersect~$\R_+$.

The vector field $1\delZ$ is well-defined on $M$. Let $f$ be obtained from
$h_*\left(1\delZ\right) = f\delz$. In a chart of a marked point
$Z_{(m,\V_m^\ell)}$, the vector field $1\delZ$ becomes a simple saddle node, so
that $f$ can be extended on $h(Z_{(m,\V_m^\ell)})$ as a meromorphic function
with a simple pole. In a chart of $\infty$, $1\delZ$ becomes a parabolic point
of codimension $k+2$ (multiplicity $k+3$). Therefore, we can extend $f$ at
infinity by $0$. Since $f$ is a meromorphic function of $\hat\Cc$ with a single
zero of order $k+3$ at infinity, it must be of the form $f = {1\over P}$, where
$P$ is a polynomial of degree $k+1$. We can centre $P$ with a translation, and
it must be monic since in the coordinate $w = {1\over z}$, ${1\over P}\delz$ has the form
$-w^{k+3}(1 + o(1)) {\del\over\del w}$ for $|w|$ small, since
$(\psi_\infty)_*\bigl(1\delZ\bigr) = -\zeta^{k+3} (1 + o(\zeta)){\del\over\del\zeta}$
and $h\circ\psi_\infty\inv(\zeta) = \zeta + o(\zeta)$.

To realize $\ol P\delz$, we consider the vector field $|h'|^{-2}\delZ$ on $M$. A
simple calculation shows that $h_*(|h'|^{-2}\delZ) = {1\over \ol{h'\circ
h\inv}}\delz$, and since $f = h'\circ h\inv$, it follows that ${1\over
\ol{h'\circ h\inv}} = \ol{P}$ on $\Cc\setminus\{\text{zeros of $P$}\}$ 
and can be extended on $\Cc$.

\section{Anti-polynomial vector fields with maximal number of heteroclinic connections}
\label{sec:non stable}

Given a fixed $(k+2)$-noncrossing tree, consider the family of structurally
stable vector fields $\ol{P_\eta} \delz$ in the generic stratum parametrized by
the analytic invariant $\eta\in\Hh^k$. When $\Im \eta_j\to 0$,
the corresponding sepal zone collapses to a curve.
Indeed, in the rectifying coordinate, $\Im\eta_j$ is the
width of the horizontal strip.
If $\Re\eta_j > 0$, then this curve contains a heteroclinic
connection $\gamma$ and $\int_\gamma P = \Re\eta_j$. If
$\Re \eta_j = 0$, then the two simple singular points on the 
boundary of the sepal zone coalesce into a double
saddle point. In both cases, the initial noncrossing tree is no longer
adequate to describe the topology of the phase portrait.

The bifurcation occuring in $\Xi_{k+1}$ are of the following type.
\begin{enumerate}
    \item Bifurcation of a heteroclinic connection. Under pertubation
        in the right direction, the heteroclinic connection breaks and gives rise
        to a sepal zone.
    \item Bifurcation of a multiple saddle point. Under perturbation 
        in the right direction, a multiple point breaks in singular
        points of smaller codimensions and gives rise to at least one
        sepal zone or heteroclinic connection.
    \item Intersection of the former two types.
\end{enumerate}

Note that it is impossible to find a closed loop made of orbits, it
would necessarily contain a non-saddle singular point.

\subsection{Case of $k$ heteroclinic connections}\label{sec:arbre ternaire}
This case can only happen when every singular point is simple.
Indeed, the graph formed by the singular points as vertices and
heteroclinic connections as edges must be a tree, and a tree
with $k$ edges must have $k+1$ vertices. Since there are
$k$ heteroclinic connections, there remains exactly $4(k+1)
- 2k = 2k+4$ outgoing and incoming separatrices. Since each
marked point at infinity must attach to at least one separatrix
and there are $2k+4$ marked points, there is exactly one separatrix
attached to each marked point.

The topological invariant of this case has a beautiful simplicity: the
separatrix graph, made of the union of heteroclinic connections, (outgoing and
incoming) separatrices, and singular points, can be extended to a ternary tree
by adding marked points at infinity as vertices. By \emph{ternary tree}, we mean
an ordered rooted tree for which each vertex is either a leaf or has three
children. In our case, we have one too many edge and vertex, so we collapse the edge corresponding to the separatrix tangent to $\R_+$ at infinity into a vertex, which we define to be the root. This way,
the tree has $3k+4$ vertices and $3k+3$ edges.

\begin{figure}
    \includegraphics{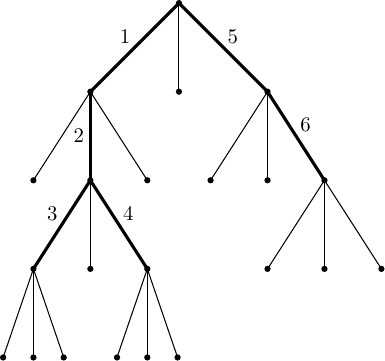}
    \caption{Ternary tree with 7 internal vertices and 6 internal edges, the latter in bold and ordered.}
    \label{fig:arbre ter ordre}
\end{figure}

We say that a vertex is an \emph{internal
vertex} if it has children. We say that an edge is an \emph{internal edge} if it
is connected to two internal vertices. We define the following total order on
the internal edges: we say that $e < f$ if either a path from the root through $f$
must first pass $e$ (see e.g.~edges $1 < 2$ on Figure~\ref{fig:arbre ter ordre}), or
a direct path from the root through $f$ is further to right than a direct
path from the root through $e$ (see e.g.~edges $1<5$ or $4<6$ on
~Figure~\ref{fig:arbre ter ordre}).

\begin{theo}\label{theo:k hetero} \begin{itemize}
    \item (Topological invariant) To every anti-polynomial vector field
        $\ol P\delz$ with $k$ heteroclinic connections, we can associate
        a unique ternary tree with $k+1$ internal vertices made with the singular points,
        the incoming and outgoing separatrices, the heteroclinic connections,
        and the marked points at infinity. The root of the tree is the separatrix tangent to $\R_+$ collapsed into a vertex.
    \item (Analytic invariants) To every anti-polynomial vector field
        $\ol P\delz$ with $k$ heteroclinic connections, we can associate
        a unique vector $\nu \in \R_+^k$, where each component is the integral
        of $P$ over an heteroclinic connection and the order is given by the 
        ternary tree.
    \item (Equivalence) Two such vector fields with the same topological invariant
        are quasi-conformally (and hence orbitally topologically) equivalent.
        Furthermore, if they also have the same analytic invariant, then
        they are analytically conjugate.
    \item (Realization) To every pair $(T_k, \nu)$, where $T_k$
        is a ternary tree with $k+1$ interior vertices and $\nu\in\R_+^k$, 
        there exists a monic centred polynomial $P$ such that
        $\ol P\delz$ has $k$ heteroclinic connections and has
        $(T_k,\nu)$ as invariants.
    \end{itemize}
\end{theo}

\begin{cor} Let $\joli H_k$ be the set of anti-polynomial vector fields
    of degree $k+1$ with $k$ heteroclinic connections. Then $\joli H_k
    /{\sim_{\text{top}}}$ has $A(k+2) = {1\over 2(k+1)+1} \binom{3(k+1)}{k+1}$ equivalence
    classes, where $\sim_{\text{top}}$ is the topological orbital equivalence fixing marked points at infinity
    of vector fields .
\end{cor}
\begin{proof}
    By Theorem~\ref{theo:k hetero}, there are as many equivalence classes of $\joli H_k
    /{\sim_{\text{top}}}$ as there are rooted ternary trees with $k+1$ internal
    vertices.
    
    In Section~\ref{subsec:counting the number of nct}, we noted that that those
    numbers are the same as the numbers we obtained for generic anti-holomorphic
    polynomial vector fields. Indeed, as discussed in that section, the number of ternary
    tree with $k$ internal edges is $A(k+2)$.
\end{proof}

\begin{proof}[Proof of Theorem~\ref{theo:k hetero}]
  The idea of the proof is the same as the Realization Theorem~\ref{theo:realisation}:
  we construct a Riemann surface isomorphic to $\hat\Cc$ with half-planes and push the vector field
  $1\delZ$ from the surface  to $\hat\Cc$. This case is even simpler, since there are
  no horizontal strips, only half-planes. The heteroclinic connections will correspond
  to intervals on the boundary of $2k$ of the half-planes.

  First, we want to embed the ternary tree in $\ol\Dd$. We start by adding a vertex
  over the root and attaching it to the root with an edge; this corresponds to the outgoing separatrix tangent to $\R_+$ at infinity that was collapsed into a vertex. Now, we order the leafs from left
  to right, i.e.~leaf 1 is the first leaf reached by a depth-first search in
  the tree by always going left, leaf 2 is the second leaf reached, and so on, and
  the added vertex over the root is leaf 0. Then we embed the tree in
  $\ol \Dd$ by mapping the $j$-th leaf on $e_j = \exp\!\left(2i\pi{j\over 2k+2}\right)$
  on $\del \Dd$, and the internal vertices to points in $\Dd$ so that the edges can
  be mapped to simple curve that do not cross each other. See Figure~\ref{fig:arbre ternaire}.

  \begin{figure}
    \subfigure[Ternary tree with added vertex atop.]{\includegraphics{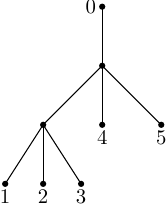}}
    \hfil
    \subfigure[Tree embedded in $\Dd$.]{\includegraphics{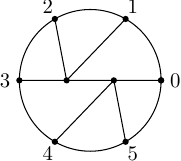}}
    \caption{Example of a tree and its embedding.}
    \label{fig:arbre ternaire}
  \end{figure}

  Now, the embedded ternary tree divides $\Dd$ in $2k+4$ simply connected domains.
  Each component is adjacent to a segment of $\del\Dd$; let $H_j^+$ be the component
  adjacent to $e_{2j-1},e_{2j}$ and $H_j^-$, the one adjacent to $e_{2j},e_{2j+1}$,
  for $1\leq j\leq k+2$. A component $H_j^\pm$ is topologically equivalent to $\pm\Hh$,
  by abuse of notation, we will use $H_j^\pm$ for the half-plane;
  a chart of the Riemann surface will be of the form $(H_j^\pm,id)$, where $id$ is the
  injection $id\colon \pm\Hh \to\Cc; z\mapsto z$. Notice that the vector field
  $1\delZ$ is well-defined on every half-plane $H_j^\pm$.

  Each $H_j^+$ (resp.~$H_j^-$) in $\Dd$ is adjacent to two edges and, say,
  $n_j^{\pm}$ internal edges. Suppose $\ell_1 < \cdots < \ell_{n_j^{\pm}}$ are the internal
  edges. Then we mark $(-\infty,0)$ with $e_{2j-1}$, $(0,\nu_{\ell_1})$ by
  $\ell_1$, $\ldots$, $(\nu_{\ell_1}+\cdots+\nu_{\ell_{n_j^{\pm}-1}},
  \nu_{\ell_1}+\ldots+\nu_{\ell_{n_j^{\pm}}})$ by $\ell_{n_j^{\pm}}$, and lastly $(\nu_{\ell_1}
  +\cdots+\nu_{\ell_{n_j^{\pm}}},\infty)$ by $e_{2j-2}$ (resp.~$e_{2j}$). We mark
  $0$ and $\nu_{\ell_1},\ldots,\nu_{\ell_{n_j^{\pm}}}$ by the corresponding vertices of the tree.

  \begin{figure}
    \includegraphics{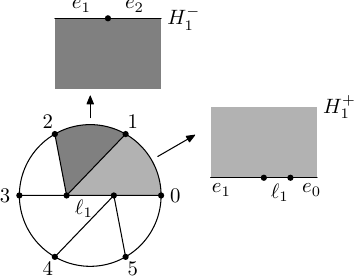}
    \caption{Example of embedded tree in $\Dd$ and the markings on half-planes.}
  \end{figure}

  We are now ready to construct the Riemann surface. Let $M$ be the Riemann surface
  obtained by sewing a marked interval of $H_j^\pm$ with the corresponding marked
  interval on another $H_m^\pm$. Equivalently, this can be seen in $\Dd$, where
  $H_j^\pm$ will share an edge with $H_m^\pm$. We can add the marked points to
  $M$ the same way we did in Section~\ref{sec:constr SR}. We also add a point at infinity
  the same way we did in that section, except that are no ends and translations.

  The end of the proof is identical to proof in Section~\ref{sec:dem realisation}:
  a well-chosen conformal mapping $h\colon \widetilde M \to \hat\Cc$ 
  will push $1\delZ$ to ${1\over P}\delz$, with $P$ a monic centred
  polynomial of degree $k+1$. Moreover, $h$ pushes $|h|^{-2}\delZ$
  to $\ol P\delz$. The phase portrait of $\ol P\delz$ must have $k$
  heteroclinic connections by construction, and its invariants must be
  $(T_k,\eta)$.
\end{proof}

\def\centrer#1{$\vcenter{\hbox{#1}}$}
\begin{figure}
    \subfigure[Bifurcation diagram of $\dot z = \ol{z^2 - \ep}$.]{\includegraphics{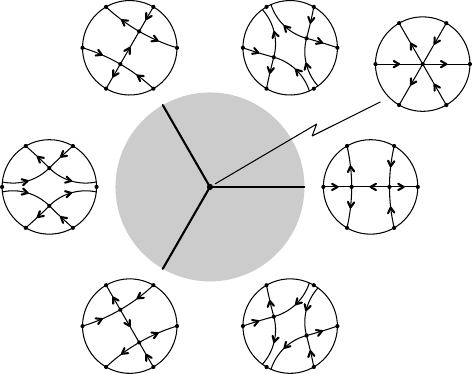}} 
    \hfil\hfil
    \subfigure[Topological invariant associated with each stratum of codimension~0 and~1.]{\includegraphics{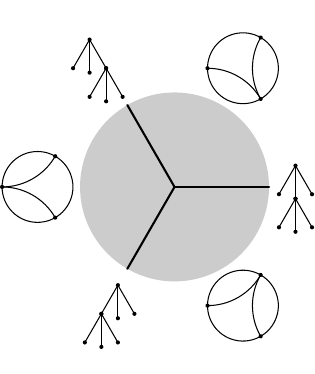}}
    \caption{Bifurcation diagram for quadratic anti-polynomial vector fields.}
\end{figure}

\subsection{Bifurcation diagram for quadratic anti-polynomial vector fields}
As a simple application of the classification of generic anti-polynomial vector
fields and anti-polynomial vectors with a maximal number of heteroclinic
connection, we can describe the complete bifurcation diagram of the family $\dot
z = \ol{z^2 - \ep} =: \ol{P_\ep(z)}$, with $\ep\in\Cc$. The only possible
non-structurally stable objects are a double saddle point, which only occur when
$\ep=0$, and a heteroclinic connection, which happens if and only if
$\int_{-\sqrt \ep}^{\sqrt\ep} P_\ep$ is real, where $\pm\sqrt\ep$ are the roots
of $P_\ep$.

A simple calculation shows that
$$
    \int_{-\sqrt \ep}^{\sqrt\ep}(z^2 - \ep)\d z
    = -{4\ep\sqrt\ep\over 3}.
$$
This is in $\R^*$ if and only if $\ep\not=0$ and $\arg\ep = {2i\pi n\over 3}$,
$n\in \mathbb Z$. 
The parameter space is divided by three rays $r$, $re^{2i\pi\over 3}$
and $re^{4i\pi\over 3}$, $r\in [0,\infty)$, and each open component outside the
rays is a generic stratum. 

Starting with $\arg\ep = 0$, the phase portrait has a heteroclinic
connection oriented from $\sqrt\ep$ to $-\sqrt\ep$ (since $\int_{-\sqrt\ep}^
{\sqrt\ep} P_\ep < 0$) on the real interval and
the phase portrait is symmetric with respect to $\R$. When
$\arg\ep$ increases in the interval $(0,\pi)$, the singular
point $\sqrt\ep$ rotates counter-clockwise around 0 until
it reaches $i\R_+$. At $\arg\ep = {2\pi\over 3}$, there
is again a heteroclinic connection oriented from $-\sqrt\ep$ to
$\sqrt\ep$ (since $\int_{-\sqrt\ep}^{\sqrt\ep} P_\ep > 0$).

When $\arg\ep$ decreases from 0 to $-\pi$, the singular
point $\sqrt\ep$ rotates clockwise around 0 until it
reaches $-i\R_+$ at $\arg\ep = -\pi$. At $\arg\ep = -{2\pi\over 3}$,
the phase portrait has a heteroclinic connection oriented
from $-\sqrt\ep$ to $\sqrt\ep$ (since $\int_{-\sqrt\ep}^{\sqrt\ep} P_\ep > 0$).

\bibliographystyle{alpha}
\bibliography{ref}

\end{document}